\documentclass[12pt,a4paper]{amsart}
\usepackage{graphics}
\usepackage{epsfig}
\usepackage{graphicx}
\usepackage[all]{xy}
\theoremstyle{plain}
\usepackage{amssymb}

\advance\hoffset-20mm \advance\textwidth40mm

\newtheorem{theorem}{Theorem}
\newtheorem{lemma}{Lemma}
\newtheorem*{theo*}{Theorem}
\newtheorem{proposition}{Proposition}
\newtheorem{corollary}{Corollary}
\theoremstyle{definition}
\newtheorem{definition}{Definition}
\newtheorem*{definition*}{Definition}
\newtheorem{example}{Example}
\newtheorem{remark}{Remark}

\def\AA{{\mathbb A}}

\def\XX{{\mathbb X}}
\def\KK{{\mathbb K}}
\def\ZZ{{\mathbb Z}}
\def\HH{{\langle H \rangle}}
\def\Ox{{\mathcal{X}}}
\def\Of{{\mathcal{O}}}
\def\Or{{\mathcal{R}}}

\def\St{\mathop{\rm Stab}}

\def\Spec{\mathop{\rm Spec}}
\def\SL{\mathop{\rm SL}}

\def\PSL{\mathop{\rm PSL}}
\def\GL{\mathop{\rm GL}}
\def\Res{\mathop{\rm Res}}
\def\Norm{\mathop{\rm Norm}}
\def\Cl{\mathop{\rm Cl}}

\begin{document}
\sloppy
\title[Torsors over Luna strata]
{Torsors over Luna strata}
\author[Ivan Arzhantsev]{Ivan V. Arzhantsev}
\thanks{Supported by RFBR grants 09-01-00648-a.}
\address{Department of Higher Algebra, Faculty of Mechanics and Mathematics,
Moscow State University, Leninskie Gory 1, GSP-1, Moscow, 119991,
Russia } \email{arjantse@mccme.ru}
\date{\today}
\begin{abstract}
Let $G$ be a reductive group and $X_H$ be a Luna stratum
on the quotient space $V/\!/G$ of a rational $G$-module $V$.
We consider torsors over $X_H$ with both non-commutative and
commutative structure groups. It allows us to describe the divisor
class group and the Cox ring of a Luna stratum under mild assumptions.
This approach gives a simple cause why many Luna strata
are singular along their boundary.
\end{abstract}
\maketitle
\section*{Introduction}

Consider a reductive group $G$ and a rational finite-dimensional $G$-module $V$
over an algebraically closed field $\KK$ of characteristic zero.
The inclusion of the algebra of invariants $\KK[V]^G$ in the polynomial algebra
$\KK[V]$ gives rise to the quotient morphism $\pi \colon V \to V/\!/G:=\Spec\,\KK[V]^G$.
In \cite{Lu73}, D.~Luna introduced a stratification of the quotient space $X:=V/\!/G$
by smooth locally closed subvarieties. Every stratum is defined by the conjugacy class
of the isotropy group $H$ of a closed $G$-orbit in the fibre of $\pi$. With any
such subgroup $H$ one associates the group $W=N_G(H)/H$.
Then the stratum $X_H$ comes with a $W$-torsor
$\pi \colon V^{\HH} \to X_H$, where $V^{\HH}$ is an open subset in the
subspace $V^H$ of $H$-fixed vectors.

The aim of this work is to use this observation to describe the Cox rings of
Luna strata and to apply them to the study of geometric properties of the strata
and their closures. It is well known that a smooth variety $Y$ with a finitely generated
divisor class group $\Cl(Y)$ admits a canonical presentation as a geometric quotient of
a quasiaffine variety by an action of a quasitorus $Q$. This presentation also is known
as a universal torsor over $Y$. The quasiaffine variety mentioned above is the
relative spectrum of the so-called Cox sheaf on $Y$ and $Q$ is a direct product of
a torus and a finite abelian group such that the group of characters $\XX(Q)$ is
identified with $\Cl(Y)$.

The $W$-torsor over a stratum $X_H$ is not a good candidate for a universal
torsor, because the group $W$ in general is not commutative. A naive idea is to
consider the commutant $S=[W,W]$, the factor group $Q=W/S$, and to decompose
the $W$-torsor into two steps, namely
$$
\xymatrix{
V^{\HH}\ar[rr]^{/W} \ar[dr]^{/S}& & X_H \\
& V^{\HH}/S \ar[ur]^{/Q} &
}
$$
Note that the groups $H$, $W$, $S$, and $Q$ are reductive.
In Theorem~\ref{main} we show that this indeed gives
a universal $Q$-torsor $V^{\HH}/S \to X_H$ under some
"codimension two" assumption on the stratum. In particular, the divisor
class group of $X_H$ is isomorphic to $\XX(Q)$.

It turns out that in many cases a Luna stratum is singular
along its boundary, i.e. the stratum coincides with the smooth locus
of its closure in $X$, see \cite[Theorem~8.1]{Ri}, \cite[Theorem~7]{LBP}
and \cite[Theorem~1.2]{KR}. In \cite{KR}, J.~Kuttler and Z.~Reichstein use
this property to prove
that quite often the Luna stratification is intrinsic. The latter means that
every automorphism of $X$ as of an abstract affine variety preserves the stratification
possibly permuting the strata. In our terms the reason for $X_H$ to be
singular along its boundary may be formulated as follows. Since the canonical
quotient presentation of the normalization $\Norm(\overline{X_H})$
of the closure of the stratum is given in Theorem~\ref{main} by the quotient morphism
$p \colon V^H/\!/S \to V^H/\!/W$ with the acting group $Q$, it is natural to expect
that $V^{\HH}/S$ is the preimage of the principal Luna stratum in $V^H/\!/W$ for
the $Q$-variety $V^H/\!/S$. Thus the fibre of $p$ over a point $x\in V^H/\!/W$
is a free $Q$-orbit if and only if $x$ is in~$X_H$, and this looks like a characterization
of smooth points on $V^H/\!/W$. We realize this approach in Theorem~\ref{T2}
under the assumption that the group $W$ is commutative.

The text is organized as follows. In Sections~\ref{sec1} and \ref{sec2}
we recall basic facts on the Luna stratification and on Cox rings respectively.
The class of admissible Luna strata which we are going to deal with in
the theorems is defined and discussed in Section~\ref{sec3}. In particular,
we prove in Proposition~\ref{Prwide} that if $H$ is the isotropy group of
a point with a closed $G$-orbit in a module $V$, then the Luna stratum
corresponding to $H$ and to the $G$-module $V^{\oplus k}$ with $k\ge 2$
is admissible. The idea to obtain strata with good properties by replacing
a module $V$ by $V^{\oplus k}$ is taken from \cite{KR}. In Section~\ref{sec4}
we formulate the main result (Theorem~\ref{main}) and obtain some corollaries.
The proof of Theorem~\ref{main} is given in Section~\ref{sec5}. Finally,
Section~\ref{sec6} is devoted to examples.


\section{The Luna stratification}
\label{sec1}

In this section we recall basic facts on the Luna stratification
obtained in~\cite{Lu73}, see also \cite[Section~6]{PV} and \cite{Sch}.
Let $G$ be a reductive affine algebraic group over
an algebraically closed field $\KK$ of characteristic zero and
$V$ be a rational finite-dimensional $G$-module. Denote
by $\KK[V]$ the algebra of polynomial functions on $V$ and
by $\KK[V]^G$ the subalgebra of $G$-invariants. Let $V/\!/G$
be the spectrum of the algebra $\KK[V]^G$. The inclusion
$\KK[V]^G\subseteq\KK[V]$ gives rise to a morphism $\pi \colon
V\to V/\!/G$ called the {\it quotient morphism} for the $G$-module $V$.
It is well known that the morphism $\pi$ is a categorical quotient
for the action of the group $G$ on $V$ in the category of
algebraic varieties, see~\cite[4.6]{PV}. In particular, $\pi$ is surjective.

The affine variety $X:=V/\!/G$ is irreducible and normal.
It is smooth if and only if the point $\pi(0)$ is smooth
on $X$. In the latter case $X$ is an affine space.
Every fibre of the morphism $\pi$ contains a unique closed $G$-orbit.
For any closed $G$-invariant subset $A\subseteq V$ its image $\pi(A)$
is closed in $X$. These and other properties of the quotient morphism
may be found in~\cite[4.6]{PV}.

By Matsushima's criterion, if an orbit $G\cdot v$ is closed
in $V$, then the isotropy group $\St(v)$ is reductive,
see \cite{Ma}, \cite{On}, or~\cite[4.7]{PV}. Moreover,
there exists a finite collection $\{H_1,\ldots,H_r\}$
of reductive subgroups in $G$ such that if an orbit $G\cdot v$
is closed in $V$, then $\St(v)$ is conjugate to one of these
subgroups. This implies that every fibre of
the morphism $\pi$ contains a point whose isotropy group coincides
with some $H_i$.

For every isotropy group $H$ of a closed $G$-orbit in $V$ the subset
$$
V_H := \{ w\in V \, ; \, \text{there exists}\,
v\in V \, \text{such that} \
\overline{G\cdot w} \supset G\cdot v =\overline{G\cdot v} \
\text{and} \, \St(v)=H\}
$$
is $G$-invariant and locally closed in $V$.
The image $X_H:=\pi(V_H)$ turns out to be a smooth
locally closed subset of $X$. In particular, $X_H$ is
a smooth quasiaffine variety.

By {\it stratification} of a variety $X$ we mean a decomposition
of $X$ into disjoint union of smooth locally closed subsets.

\begin{definition}
The stratification
$$
X \, = \, \bigsqcup_{i=1}^r \, X_{H_i}
$$
is called the {\it Luna stratification} of the quotient space $X$.
\end{definition}

There is a unique open dense stratum called the {\it principal
stratum} of $X$. A stratum $X_{H_i}$  is contained in the closure of a stratum
$X_{H_j}$ if and only if the subgroup $H_i$ contains a subgroup conjugate to $H_j$.
This induces a partial ordering of the set of strata compatible with the (reverse)
ordering on the set of conjugacy classes.

Consider the subsets
$$
V^H=\{ v\in V \, ; \, H\cdot v =v \} \quad \text{and}
\quad V^{\HH}=\{v\in V \, ; \, \St(v)=H \, \text{and} \,
G\cdot v = \overline{G\cdot v}\}.
$$
Then $V^{\HH}$ is an open subset of $V^H$. Moreover, the restriction
of the morphism $\pi$ to $V^H$ maps $V^{\HH}$ to $X_H$ surjectively.

Let $N_G(H)$ be the normalizer of the subgroup $H$ in $G$.
Since $H$ is reductive, the subgroup $N_G(H)$ is reductive and
the connected component $N_G(H)^0$ is the product of $H^0$
and $C_G(H)^0$, where $C_G(H)^0$ is the connected component of the centralizer
$C_G(H)$ of $H$ in $G$,
see \cite[Lemma~1.1]{LR}. The group $N_G(H)$ preserves the subspace $V^H$.
Moreover, the kernel of the $N_G(H)$-action on $V^H$ is $H$, and
we get an effective action of the (reductive) group $W := N_G(H)/H$
on $V^H$. By \cite{Lu75}, for any point $v\in V^H$ the orbit
$G\cdot v$ is closed in $V$ if and only if the orbit $W\cdot v$
is closed in $V^H$. In particular, $V^{\HH}$ is the union of
closed free $W$-orbits in~$V^H$.

\smallskip

The following definition first appeared in \cite{Po} plays an important
role in Invariant Theory.

\begin{definition}
An action of a reductive group $F$ on an affine variety $Z$ is
{\it stable} if there exists an open dense subset $U\subseteq Z$ such that
the orbit $F\cdot z$ is closed in $Z$ for any $z\in U$.
\end{definition}

Equivalently, an $F$-action on $Z$ is stable if the general fibre of the
quotient morphism $\pi \colon Z \to Z/\!/F:=\Spec\,\KK[Z]^F$
is an $F$-orbit, i.e. there is an open dense subset $U' \subseteq Z/\!/F$
such that $\pi^{-1}(u)$ is an $F$-orbit for any $u\in U'$.

\smallskip

Since $V^{\HH}$ is open in $V^H$, the action of $W$ on $V^H$ is stable.
This action enjoys an additional nice property, namely, the isotropy group of
the general point in $V^H$ is trivial.

Summarizing, we observe that the
restriction of the quotient morphism $\pi$ to $V^{\HH}$ defines
a $W$-torsor $\pi \colon V^{\HH} \to X_H$.

Restricting the algebra of invariants $\KK[V]^G$ to the subspace
$V^H$, we get a subalgebra $\Res_H(\KK[V]^G)$ in $\KK[V^H]^W$.
The closure $\overline{X_H}$ may be identified with $\Spec\,\Res_H(\KK[V]^G)$.
By~\cite{Lu75}, the morphism
$$
\pi_{H} \colon V^H/\!/W \to \overline{X_H}
$$
given by the inclusion $\Res_H(\KK[V]^G) \subseteq \KK[V^H]^W$
is the normalization morphism. For the principal stratum,
the restriction $\Res_H$ is injective, the closure $\overline{X_H}=X$
is normal, and we get an isomorphism $\KK[V]^G \to \KK[V^H]^W$.
Several results on normality of the closures of Luna strata and on strata
of codimension one may be found in \cite{Sch} and \cite{Shm}.


\section{Cox rings and characteristic spaces}
\label{sec2}

Let $X$ be a normal algebraic variety with finitely generated divisor class group $\Cl(X)$.
Assume that any regular invertible function $f\in \KK[X]^{\times}$ is constant.
Roughly speaking, the {\it Cox ring} of $X$ may be defined as
$$
R(X) \, := \, \bigoplus_{D\in \Cl(X)} \Gamma(X, \Of_X(D)).
$$
In order to obtain a multiplicative structure on $R(X)$ some technical work is needed,
especially when the group $\Cl(X)$ has torsion. We refer for details to
\cite[Section~4]{ADHL} and \cite{Ha}. In a similar way one defines the {\it Cox sheaf} $\Or$,
which is a sheaf of $\Cl(X)$-graded algebras on $X$. The Cox ring $R(X)$ is the ring of global
sections of the sheaf $\Or$.
Assume that this sheaf is locally of finite type.
Then the relative spectrum $\widehat{X} := \Spec_X(\Or)$ is a quasiaffine variety.

By definition, a {\it quasitorus}~\footnote{Also called a diagonalizable group.} is a commutative reductive
algebraic group. It is not difficult to check that any quasitorus is a direct product of a torus and
a finite abelian group. Since  $\Or$ is a sheaf of $\Cl(X)$-graded algebras, the variety $\widehat{X}$ comes with an action
of  a quasitorus $Q_X$ whose group of characters $\XX(Q_X)$ is identified with $\Cl(X)$.

Now we need a notion from Geometric Invariant Theory. Let $G$ be a reductive group and $Z$ be a normal $G$-variety.
An invariant morphism $\pi \colon Z \to X$ is a {\it good quotient}, if $\pi$ is affine and the pullback map
$\pi^{*} \colon \Of_X \to \Of_Z^G$ is an isomorphism. An example of a good quotient is the quotient morphism
$\pi \colon V\to V/\!/G$ discussed in Section~\ref{sec1}. Moreover, for any open subset $U\subseteq V/\!/G$
the restriction $\pi \colon \pi^{-1}(U) \to U$ is a good quotient as well.

It turns out that the $Q_X$-action on $\widehat{X}$ admits a good quotient $p_X \colon \widehat{X} \to X$.
This morphism is called the {\it characteristic space} over $X$. If the variety $X$ is smooth,
then the $Q_X$-action on $\widehat{X}$ is free, and the map $p_X \colon \widehat{X} \to X$ is
a $Q_X$-torsor. Moreover, this is a {\it universal torsor} over $X$ in the sense of~\cite{Sk}.

\smallskip

The following definition appeared in \cite[Definition~6.4.1]{ADHL}, see also \cite[Definition~1.37]{Ha}.

\begin{definition} \label{defstst}
An action of a reductive group $F$ on an affine variety $Z$ is said to be
{\it strongly stable} if there exists an open dense invariant subset $U\subseteq Z$ such that
\begin{enumerate}
\item
the complement $Z\setminus U$ is of codimension at least two in $Z$;
\item
the group $F$ acts freely on $U$;
\item
for every $z\in U$ the orbit $F\cdot z$ is closed in $Z$.
\end{enumerate}
\end{definition}

One may check that the $Q_X$-action on $\widehat{X}$ defined above is strongly stable
with $U$ being the preimage $p_X^{-1}(X_{\text{reg}})$ of the smooth locus of $X$.

\begin{definition}
A normal variety $Y$ equipped with an action of a quasitorus $Q$ is
called {\it $Q$-factorial} if any $Q$-invariant Weil divisor on $Y$
is principal.
\end{definition}

Again the characteristic space $\widehat{X}$ is a $Q_X$-factorial variety,
see \cite[Lemma~5.3.6]{ADHL} and \cite[Theorem~1.14]{Ha}.

\smallskip

The following theorem given in \cite[Theorem~6.4.3]{ADHL}
and \cite[Theorem~1.40]{Ha} shows that strong stability and
$Q$-factoriality define exactly the class of characteristic spaces.
This result is our main tool in the proof of Theorem~\ref{main}.

\begin{theorem} \label{Tquot}
Let a quasitorus $Q$ act on a normal quasiaffine variety $\Ox$ with
a good quotient $p \colon \Ox \to X$. Assume that $\KK[\Ox]^{\times}=\KK^{\times}$
holds, $\Ox$ is $Q$-factorial and the $Q$-action is strongly stable. Then there is
a commutative diagram
$$
\xymatrix{
\Ox \ar[rr]^{\mu}_{\cong} \ar[dr]_{p} & & \widehat{X} \ar[dl]^{p_X} \\
& X &
}
$$
where the quotient space $X$ is a normal variety with $\KK[X]^{\times}=\KK^{\times}$,
we have $\Cl(X)=\XX(Q)$, the morphism $p_X \colon \widehat{X} \to X$ is a characteristic
space over $X$, and the isomorphism $\mu \colon \Ox \to \widehat{X}$ is equivariant with
respect to actions of $Q=Q_X$.
\end{theorem}


\section{Admissible Luna strata}
\label{sec3}

Let us fix settings and notation for the remaining part of the text.
Let $G$ be a reductive group, $V$ be a rational finite-dimensional $G$-module,
$\pi \colon V \to V/\!/G:=\Spec\,\KK[V]^G$ be the quotient morphism with the
quotient space $X:=V/\!/G$, and $H$ be an isotropy group of a closed $G$-orbit in $V$.
Define the group $W=N_G(H)/H$, its commutant $S:=[W,W]$ and the factor group $Q:=W/S$.
Recall that the groups $W$ and $S$ are reductive. Moreover, $S$ is semisimple
provided $W$ is connected.
Being reductive and commutative, the group $Q$ is a quasitorus.

\smallskip

We adapt Definition~\ref{defstst} to the situation we are dealing with.

\begin{definition}
A triple $(G,V,H)$
is {\it admissible}, if the action of the group $W$ on $V^H$ is
strongly stable.  An {\it admissible} stratum is the Luna stratum $X_H$
corresponding to an admissible triple $(G,V,H)$.
\end{definition}

\begin{remark} \label{rem1}
A triple $(G,V,H)$ is admissible if and only if the complement
of $V^{\HH}$ in $V^H$ has codimension at least two. Indeed,
we know from Section~\ref{sec1} that the action of $W$ on $V^{\HH}$
is free, all $W$-orbits on $V^{\HH}$ are closed in $V^H$ and the subset
$V^{\HH}$ is characterized by these properties. In other words, the subset
$V^{\HH}$ is the preimage of the principal Luna stratum in the $W$-module $V^H$.
This shows that in order to check admissibility it suffices to work
with principal strata.
\end{remark}

Since $\overline{X_H}=\pi(V^H)$, the complement of an admissible stratum
in its closure has codimension at least two. In this case the algebra of regular
functions $\KK[X_H]$ coincides with $\KK[\Norm(\overline{X_H})]$,
and the latter algebra is isomorphic to $\KK[V^H]^W$, see Section~\ref{sec1}.
In particular, $\KK[X_H]^{\times}=\KK^{\times}$ holds for any admissible stratum.

Following \cite{KR} we say that a Luna stratum $X_H$ is {\it singular along its boundary} if
the singular locus of the closure $\overline{X_H}$ is precisely the complement
$\overline{X_H}\setminus X_H$. This property is known for the principal
stratum on the quotient variety of the space of representations of a quiver
except for some low dimensional anomalies \cite[Theorem~7]{LBP}. Moreover,
there are many modules where every Luna stratum is singular along its
boundary, see~\cite[Theorem~1.2]{KR}. In particular, this is the case for
the module $V^{\oplus r}$, where $V$ is any rational $G$-module and
$r\ge 2\dim(V)$, or where $V$ is the adjoint module and $r\ge 3$. If the closure
of a Luna stratum is normal, then "singular along its boundary" implies that
the boundary has codimension at least two. This implication always holds
for the principal stratum. If we know additionally that the morphism
$\pi \colon V \to X$ does not contract invariant
divisors~\footnote{This is the case if $G$ is semisimple.},
then the principal stratum singular along its boundary is admissible.

\smallskip

The following proposition provides us with a wide class of admissible strata.

\begin{proposition} \label{Prwide}
Let $V$ be a rational $G$-module and $H$ be an isotropy group of a closed $G$-orbit in $V$.
For any positive integer $k$ consider the module
$V^{\oplus k}=V\oplus\cdots\oplus V$ ($k$ times) with the diagonal $G$-action.
Then the subgroup $H$ is an isotropy group of a closed $G$-orbit in $V^{\oplus k}$ and for $k\ge 2$
the triple $(G, V^{\oplus k}, H)$ is admissible.
\end{proposition}

\begin{proof}
Since $H$ is an isotropy group of a closed $G$-orbit in $V$, there is a point $v\in V$ such that
the orbit $G\cdot v$ is closed and the isotropy group $\St(v)$ coincides with $H$.
The $G$-orbit of the point $(v,\ldots,v)\in V^{\oplus k}$ is contained
in the diagonal $\Delta V \subseteq V^{\oplus k}$. Since $\Delta V \cong V$, this
orbit is closed in $\Delta V$ and thus in $V^{\oplus k}$. This proves that $H$ is
an isotropy group of a closed $G$-orbit in $V^{\oplus k}$.

Consider the diagonal action of $W$ on $(V^{\oplus k})^H=(V^H)^{\oplus k}$
with $k\ge 2$. We claim that if at least one component $v_i$ of a point
$v=(v_1,\ldots,v_k)$ is contained in $V^{\HH}$, then $v$ is contained in
$(V^{\oplus k})^{\HH}$. Indeed, assume that $v_1\in V^{\HH}$. Then the $W$-isotropy group
of $v_1$, and thus of $v$, is trivial. The orbit $W\cdot v$ is contained
in a closed subset $(W\cdot v_1) \times (V^H)^{\oplus (k-1)}$. If this orbit
is not closed, its closure contains a $W$-orbit of smaller dimension.
On the other hand, any point in $(W\cdot v_1) \times (V^H)^{\oplus (k-1)}$
has the trivial isotropy group.

We conclude that the complement of $(V^{\oplus k})^{\HH}$ in $(V^{\oplus k})^H$
is contained in the set of points $(v_1,\ldots,v_k)$, where every $v_i$ is not in $V^{\HH}$.
Clearly, the latter set has codimension at least two in $(V^H)^{\oplus k}$.
\end{proof}

Finally let us consider the case when the group $G$ is finite. By Remark~\ref{rem1}
it suffices to know when the principal stratum is admissible. Recall that a linear operator
is called {\it pseudoreflection} if it has a finite order and its subspace of fixed
vectors is a hyperplane.

\begin{proposition} \label{profin}
Let $G$ be a finite group and $V$ be a $G$-module. Then the principal Luna stratum
in $X$ is admissible if and only if the image of $G$ in $\GL(V)$
does not contain pseudoreflections.
\end{proposition}

\begin{proof}
Any $G$-orbit in $V$ is closed, so the principal stratum is admissible
if and only if the set of points in $V$ with non-trivial isotropy groups has
codimension at least two. This means that for every $g\in G$, $g\ne e$, its subspace of
fixed vectors has codimension at least two, or, equivalently,
$g$ is not a pseudoreflection.
\end{proof}


\section{The Cox ring of an admissible stratum}
\label{sec4}

We preserve the settings of the previous section.
The quasitorus $Q$ acts on the quotient space $V^H/\!/S$  and this action defines
an $\XX(Q)$-grading on the algebra $\KK[V^H]^S$.

\begin{theorem} \label{main}
Let $G$ be a reductive group, $V$ be a rational finite-dimensional $G$-module,
and $H$ be an isotropy group of a closed $G$-orbit in $V$.
Assume that the stratum $X_H$ is admissible.
Then the $Q$-quotient morphism
$$
p \, \colon \, V^H/\!/S \longrightarrow V^H/\!/W \cong \Norm(\overline{X_H})
$$
is a characteristic space over the normalization $\Norm(\overline{X_H})$.
Moreover, if $q \colon V^H \to V^H/\!/S$ is the $S$-quotient morphism and $Y:=q(V^{\HH})$, then
the restriction of $p$ to $Y$ defines a $Q$-torsor
$$
p \, \colon Y \, \longrightarrow X_H,
$$
which is a universal torsor over $X_H$.
In particular, the divisor class group $\Cl(X_H)$ may be identified with $\XX(Q)$,
and the Cox ring $R(X_H)$ is isomorphic to $\KK[V^H]^S$ as an $\XX(Q)$-graded ring.
\end{theorem}

\begin{corollary}
The Cox ring of an admissible Luna stratum is finitely generated.
\end{corollary}

\begin{proof}
Since the group $S$ is reductive, the algebra of invariants $\KK[V^H]^S$
is finitely generated.
\end{proof}

\begin{corollary} \label{cor2}
If a Luna stratum $X_H$ is admissible and the group $W$ is connected, then
the group $\Cl(X_H)$ is free and the Cox ring $R(X_H)$ is factorial.
\end{corollary}

\begin{proof}
Since the group $W$ is connected and reductive, the factor group $Q$ is a torus,
and the group $\XX(Q)\cong\Cl(X_H)$ is a lattice. Moreover, the commutant $S$ is semisimple
and the algebra of invariants $\KK[V^H]^S\cong R(X_H)$ is factorial, see \cite[Theorem~3.17]{PV}.
\end{proof}

\begin{remark}
In more general settings, the Cox ring of a variety with a free divisor class group is
a unique factorization domain, see. e.g. \cite[Theorem~1.14]{Ha}.
\end{remark}

\begin{corollary}
Let $(G,V,H)$ be an admissible triple. Then the divisor class group of the variety
$\Norm(\overline{X_H})$ is isomorphic to $\XX(Q)$, while the Picard group of
$\Norm(\overline{X_H})$ is trivial.
\end{corollary}

\begin{proof}
The first statement follows directly from Theorem~\ref{main}.
For the second one we note that the origin in $V^H$ is a $W$-fixed
point which projects to a $Q$-fixed point on $V^H/\!/S$. Since for
the characteristic space $p \colon V^H/\!/S \to \Norm(\overline{X_H})$
we have a $Q$-fixed point above, it follows from \cite[Corollary~1.36]{Ha}
that the Picard group of $\Norm(\overline{X_H})$ is trivial.
\end{proof}

The following theorem shows that many admissible strata are singular along
their boundary.

\begin{theorem} \label{T2}
Let $G$ be a reductive group, $V$ be a rational finite-dimensional $G$-module, and
$H$ be an isotropy group of a closed $G$-orbit in $V$. Assume that the stratum $X_H$ is admissible
and the group $W$ is commutative. Then the stratum $X_H$ is singular along
its boundary.
\end{theorem}

\begin{proof}
Under our assumptions the group $S$ is trivial and so $Q$ coincides with $W$.
By Theorem~\ref{main}, the quotient morphism $p \colon V^H \to V^H/\!/W$ is the
characteristic space over $V^H/\!/W \cong \Norm(\overline{X_H})$. Since the
variety $V^H$ is smooth, a point $x\in \Norm(\overline{X_H})$ is smooth if and only
if the fibre $p^{-1}(x)$ consists of a unique free $W$-orbit,
see~\cite[Proposition~4.12]{Ha-II}. This shows that the preimage of the smooth locus
in $\Norm(\overline{X_H})$ coincides with $V^{\HH}$ and so the smooth locus
of $\Norm(\overline{X_H})$ is $X_H$. This implies that the smooth locus
of $\overline{X_H}$ is $X_H$ too.
\end{proof}

Finally let us specialize Theorem~\ref{main} to the case when the group $G$ is a torus $T$.
Any rational $T$-module $V$ admits the weight decomposition
$$
V=\bigoplus_{\lambda\in\XX(T)} V(\lambda), \quad \text{where} \quad
V(\lambda)=\{v\in V \, ; \, t\cdot v =\lambda(t)v\}.
$$
The {\it multiplicity} of a weight $\lambda$ in $V$ is defined as $\dim V(\lambda)$.

\begin{proposition} \label{PrPr}
Let $V$ be a rational $T$-module such that every weight in $V$ has multiplicity at least two.
Then every Luna stratum $X_H$ is admissible, the divisor class group of $X_H$ is a subgroup
of the lattice $\XX(T)$ and the Cox ring $R(X_H)$ is polynomial.
\end{proposition}

\begin{proof}
For any subgroup $H$ in $T$ which is an isotropy group of a closed $G$-orbit in $V$ the group $W$ is just $T/H$.
Thus it is commutative, $S=\{e\}$, and $Q=W$. Since
$V^H=\bigoplus V(\lambda)$ with $\lambda|_H=1$, the multiplicity of any weight of
the torus $W$ in $V^H$ is at least two.  To get admissibility, it suffices to show
that an effective stable linear $W$-action satisfying this condition is strongly stable.
For any vector $v\in V^H$ define its weight system as the set of $\lambda$'s such that
the component of $v$ in $V(\lambda)$ is non-zero. Then the isotropy group of $v$ in $W$ is trivial
if and only if the weight system generates the lattice $\XX(W)$. Moreover, the orbit
$W\cdot v$ is closed in $V^H$ if and only if the convex hull of the weight system contains
zero in its relative interior. By assumptions these properties are fulfilled if all weight
components of $v$ are non-zero. As all weight multiplicities are at least two,
the codimension of the complement
of $V^{\HH}$ in $V^H$ is at least two as well, and the stratum $X_H$ is admissible.
Further,
$$
\Cl(X_H) \, \cong \, \XX(Q) \, = \, \XX(W) \, = \, \XX(T/H),
$$
and $\XX(T/H)$ is a subgroup of $\XX(T)$. Finally,
$$
R(X_H) \, \cong \, \KK[V^H]^S \, = \, \KK[V^H]
$$
is a polynomial algebra.
\end{proof}

\begin{remark}
More generally, every Luna stratum as in Theorem~\ref{T2}
is a toric variety, so its Cox ring is polynomial~\cite{Cox}.
The rays of the polyhedral cone defining the affine
toric variety $V^H/\!/W$ are obtained via Gale duality from
the weights of the quasitorus $W$ on $V^H$,
see~\cite[Chapter II, 2.1]{ADHL}, while the stratum $X_H$
may be recovered as the smooth locus of $V^H/\!/W$.
\end{remark}

\begin{corollary}
Under the assumptions of Proposition~\ref{PrPr}, every Luna stratum in
$V/\!/T$ is singular along its boundary. In particular, the Luna stratification
in $V/\!/T$ is intrinsic.
\end{corollary}

\begin{proof}
The first statement follows from Theorem~\ref{T2} and Proposition~\ref{PrPr}.
For the second one, see~\cite[Lemma~3.1]{KR}.
\end{proof}


\section{Proof of Theorem~\ref{main}}
\label{sec5}

We get the result by checking all the conditions of Theorem~\ref{Tquot}.

\begin{proposition}
Let $(G,V,H)$ be an admissible triple. Then the action of the quasitorus
$Q$ on $V^H/\!/S$ is strongly stable.
\end{proposition}

\begin{proof}
We begin with the following general observation.

\begin{lemma} \label{Lemstr}
Let $F$ be a reductive group, $S=[F,F]$ be its commutant and $Q$ be the
factor group $F/S$. If an action of the group $F$ on an irreducible affine
variety $Z$ is strongly stable, then the induced action of $Q$ on $Z/\!/S$
is strongly stable as well.
\end{lemma}

\begin{proof}
Let $U$ be an open subset of $Z$ as in Definition~\ref{defstst}.
Then the group $S$ acts on $U$ freely and all $S$-orbits on $U$
are closed in $Z$. Let $q \colon Z \to Z/\!/S$ be the quotient morphism
and $Y=q(U)$. Then $Y$ is open in $Z/\!/S$, $U=q^{-1}(Y)$ and
since $q$ is surjective, the complement of $Y$ in $Z/\!/S$
has codimension at least two. Moreover, for any $z\in U$
the image $q(F\cdot z)$ is a closed $Q$-orbit of the point $q(z)$
with trivial isotropy group. This proves that the group $Q$ acts on $Y$
freely and all $Q$-orbits on $Y$ are closed in $Z/\!/S$.
\end{proof}

Let us return to the proof of the proposition. By definition of an addmisible triple
the action of the group $W$ on $V^H$ is strongly stable, and Lemma~\ref{Lemstr}
implies that the induced action of $Q$ on $V^H/\!/S$ is strongly stable as well.
\end{proof}

\begin{proposition}
Let $(G,V,H)$ be an admissible triple.
Then the $Q$-variety $V^H/\!/S$ is $Q$-factorial.
\end{proposition}

\begin{proof}
Take any $Q$-invariant Weil divisor $D$ on $V^H/\!/S$. Without loss of generality
we may assume that $D$ is effective.
The preimage $q^*(D)$ of $D$ under the quotient morphism $q \colon V^H \to V^H/\!/S$
is a $W$-invariant divisor on $V^H$. Since $V^H$ is an affine space, any such divisor
is a principal divisor of some regular $W$-semiinvariant function. Such a function
is $S$-invariant, thus it may be considered as a regular function on $V^H/\!/S$
having $D$ as its principal divisor.
\end{proof}

Since $\KK[\Norm(\overline{X_H})]^{\times}\cong\KK[V^H/\!/W]^{\times}=\KK^{\times}$, we get the
first statement of Theorem~\ref{main} from Theorem~\ref{Tquot}. The second one follows from
the fact that the complement of $X_H$ in $\Norm(\overline{X_H})$ has codimension at least two
and \cite[Lemma~5.1.2]{ADHL}.
The proof of Theorem~\ref{main} is completed.


\section{Examples}
\label{sec6}

In this section we consider several examples. They illustrate Theorem~\ref{main}
and show that the condition for a stratum to be admissible is essential.

\begin{example}
Let $G=\SL(2)$ and $V$ be the tangent algebra of $G$ considered as adjoint $G$-module. This module
has conjugacy classes of subgroups which are isotropy groups of closed $G$-orbits,
namely the maximal torus $T$ and the group $\SL(2)$ itself.
Since $V^T$ is an affine line and the complement of $V^{\langle T\rangle}$ in $V^T$ is a point,
the triple $(\SL(2),V,T)$ is not admissible. In this case the quotient space
$X:=V/\!/G$ is an affine line $\AA^1$, and the stratum $X_T$ coincides with
$\AA^1\setminus \{0\}$. The group $W=N_G(T)/T$ has order two, so $S=\{e\}$ and
$Q=W$. The divisor class group $\Cl(\AA^1\setminus \{0\})$ is trivial, thus does not coincide
with $\XX(Q)$, while $\KK[\AA^1\setminus \{0\}]^{\times}\ne\KK^{\times}$ and the Cox ring of
$X_T$ is not defined.

Further, the principal stratum for $V\oplus V$ corresponds to $H=\{\pm E\}$. One easily checks
that it is also not admissible.
In this case $W=\PSL(2)=S$, $Q$ is trivial, and
$$
(V\oplus V)^H/\!/S=(V\oplus V)/\!/G \cong \AA^3,
$$
so all conclusions of Theorem~\ref{main} hold.
\end{example}

\begin{example} \label{ex2}
Let us take an arbitrary connected semisimple group $G$ and the adjoint $G$-module $V$.
The principal stratum for $V$ corresponds to $H=T$, where $T$ is a maximal torus in $G$,
and $V^H$ is a Cartan subalgebra. In this case $W$ is the (finite)
Weyl group associated with the group $G$ and $Q$ is its maximal commutative factor.
The principal stratum is not admissible, but starting from $k=2$ the stratum for
$V^{\oplus k}$ associated with $H=T$ is admissible. Its divisor class group does
not depend on $k$ and is isomorphic to $\XX(Q)\cong Q$. In order to describe the Cox ring
one should calculate the algebra of invariants of the diagonal action of the group
$S=[W,W]$ on the direct sum of Cartan subalgebras.
It is important to note that by~\cite[Theorem~8.1]{Ri} the principal stratum
for $V^{\oplus k}$ is singular along its boundary for any $k\ge 2$ provided
$G$ has no simple factors of rank~$1$.
\end{example}

\begin{example}
Take $G=\SL(n)$ and $V=\KK^n\oplus(\KK^n)^*$, where $\KK^n$ is the tautological
$\SL(n)$-module and $(\KK^n)^*$ is its dual. The principal stratum for $V$
corresponds to $H\cong\SL(n-1)$, and $W=Q$ is a one-dimensional torus.
Since $V/\!/G\cong \AA^1$, the principal stratum is not admissible.
But starting from $k=2$ the subgroup $H$ defines an admissible statum
for $V^{\oplus k}$ whose divisor class group is isomorphic to $\ZZ$.
By Theorem~\ref{T2}, this stratum is singular along its boundary.
\end{example}

Let us show that the Cox ring of an admissible Luna stratum is not always factorial,
cf.~Corollary~\ref{cor2}. The example below is taken from \cite[Section~3]{AG},
where the Cox ring of the quotient space of an arbitrary finite linear group is
described.

\begin{example}
Let $V=\KK^2$ and $G$ be the quaternion group
$$
Q_8=\left\{\pm E, \pm
\begin{pmatrix}
i & 0 \\
0 & -i
\end{pmatrix},
\pm
\begin{pmatrix}
0 & 1 \\
-1 & 0
\end{pmatrix},
\pm
\begin{pmatrix}
0 & i \\
i & 0
\end{pmatrix}\right\},
$$
where $i^2=-1$. By Proposition~\ref{profin}, the principal Luna stratum in $X$
corresponding to $H=\{e\}$ is admissible. As $S=[Q_8,Q_8]=\{\pm E\}$, the Cox ring
of the principal stratum $\KK[V]^S$ is the ring of functions on the two-dimensional quadratic cone.
This ring is not factorial.
\end{example}

Note that Example~\ref{ex2} also leads to non-factorial Cox rings. Say, this is already the case
for $G=\SL(3)$ and $k=2$. But the description of the Cox ring here is more complicated.


\section*{Acknowledgement}

The preparation of this text was inspired by the workshop
"Torsors: theory and applications", January 10--14, 2011,
International Centre for Mathematical Sciences, Edinburgh,
organized by Victor Batyrev and Alexei Skorobogatov.
The final version was written
during a stay of the author at the Institut Fourier, Grenoble.
He wishes to thank these institutions for generous
support and hospitality. Special thanks are due to the referee
for the careful reading, useful remarks and suggestions.


%

\begin{thebibliography}{}
%
\bibitem{ADHL}
I.V.~Arzhantsev, U.~Derenthal, J.~Hausen, and A.~Laface.
Cox rings. arXiv:1003.4229.
%
\bibitem{AG}
I.V.~Arzhantsev and S.A.~Gaifullin.
Cox rings, semigroups and automorphisms of affine algebraic varieties.
Mat. Sbornik 201 (2010), no.~1,  3--24 (Russian); English transl.: Sbornik: Math.
201 (2010), no.~1, 1--21.
%
%
\bibitem{Cox}
D.A.~Cox. The homogeneous coordinate ring of a toric variety.
J. Alg. Geom. 4 (1995), no.~1, 17--50.
%
\bibitem{Ha-II}
J.~Hausen. Cox rings and combinatorics II.
Moscow Math. J. 8 (2008), no.~4, 711--757.
%
\bibitem{Ha}
J.~Hausen. Three letures on Cox rings. This volume.
%
\bibitem{KR}
J.~Kuttler and Z.~Reichstein.
Is the Luna stratification intrinsic?
Ann. Inst. Fourier (Grenoble) 58 (2008), no.~2, 689--721.
%
\bibitem{LBP}
L. Le Bruyn and C. Procesi. Semisimple representations of quivers.
Trans. Amer. Math. Soc. 317 (1990). no.~2, 585--598.
%
\bibitem{Lu73}
D.~Luna. Slices \'etales. Bull. Soc. Math. France,
Memoire 33 (1973), 81--105.
%
\bibitem{Lu75}
D.~Luna. Adh\'erences d'orbite et invariants.
Invent. Math. 29 (1975), 231--238.
%
\bibitem{LR}
D.~Luna and R.W.~Richardson.
A generalization of the Chevalley restriction theorem.
Duke Math. J. 46 (1979), no.~3, 487--496.
%
\bibitem{Ma}
Y.~Matsushima. Espaces homog\'enes de Stein des groupes de Lie complexes.
Nagoya Math. J. 16 (1960), 205--218.
%
\bibitem{On}
A.L.~Onishchik. Complex hulls of compact homogeneous spaces.
Dokl. Akad. Nauk SSSR 130 (1960), no.~4, 726--729 (Russian);
English Transl.: Soviet Math. Dokl. 1 (1960), 88--91.
%
\bibitem{Po}
V.L. Popov. Criteria for the stability of the action of a semisimple group on a factorial manifold.
Izv. Akad. Nauk SSSR~34 (1970), 523--531 (Russian);
English transl.: Math. USSR, Izv.~4 (1970), 527--535.
%
\bibitem{PV}
V.L.~Popov and E.B.~Vinberg. Invariant Theory.
Algebraic Geometry IV, Encyclopaedia Math. Sciences, vol.~55,
Springer-Verlag Berlin, 1994, pp. 123--278.
%
\bibitem{Ri}
R.W. Richardson. Conjugacy classes of $n$-tuples in Lie algebras and algebraic groups.
Duke Math. J. 57 (1988), no.~1, 1--35.
%
\bibitem{Sch}
G.W. Schwarz. Lifting smooth homotopies of orbit spaces.
Inst. Hautes \'Etudes Sci. Publ. Math. 51 (1980), 37--135.
%
\bibitem{Shm}
D.A.~Shmel'kin. On algebras of invariants and codimension $1$
Luna strata for nonconnected groups. Geom. Dedicata 72 (1998),
no.~2, 189--215.
%
\bibitem{Sk}
A.N.~Skorobogatov. Torsors and rational points. Cambridge Tracts in Mathematics 144,
Cambridge University Press, Cambridge, 2001, 187 pp.
%
\end{thebibliography}
\end{document}